%&AMS-TeX
%\input amstex.tex
\documentstyle{amsppt}
\magnification=1200
\hsize=150truemm
\vsize=224.4truemm
\hoffset=4.8truemm
\voffset=12truemm

%\NoBlackBoxes
\NoRunningHeads

\define\C{{\bold C}}%                \C          C mathematique
 
%                \Q          Q mathematique
%                \R          R mathematique
%                \Z          Z mathematique
%                \N          N mathematique
 
\let\thm\proclaim%             \thm        \proclaim
\let\fthm\endproclaim%         \fthm       \endproclaim
 %              \inc        inclus dans      (\subset)
%\let\ds\displaystyle%          \ds         \displaystyle
%              \ev         ensemble vide   (\emptyset)

 \define\p{ P^1(\C) }%         \p      P1(C)
\define\pp{ P^2(\C) }%         \pp      P2(C)
\define\De{\Delta}%              \De     Delta
\define\Dn{\Delta_n}%              \Dn     Delta indice n
\define\dn{\delta_n}%              \dn     delta indice n
  \define\Si{\Sigma}%              \Si    Sigma
\define\Ti{\widetilde{T}}%         \Ti    T tilda
\define\DTn{\widetilde{\Delta}_n}%  \DTn    Delta indice n tilda
%%%%%%%%    numerotation des sections et enonces %%%%%%%%%%%%
\newcount\tagno
\newcount\secno
\newcount\subsecno
\newcount\stno
\global\subsecno=1
\global\tagno=0
\define\ntag{\global\advance\tagno by 1\tag{\the\tagno}}

\define\sta{\ %\uppercase\expandafter{\romannumeral
{\the\secno}.\the\stno
\global\advance\stno by 1}

\define\stas{\the\stno
\global\advance\stno by 1}

\define\sect{\global\advance\secno by 1
\global\subsecno=1\global\stno=1\
%\uppercase\expandafter{\romannumeral\the\secno}. }
{\the\secno}. }

\def\nom#1{\edef#1{{\the\secno}.\the\stno}}
\def\inom#1{\edef#1{\the\stno}}
\def\eqnom#1{\edef#1{(\the\tagno)}}

%%%%%%%%%%%%%%  macros pour la bibliographie  %%%%%%%%%%%%%%%%%%%%%%
\newcount\refno
\global\refno=0

\def\nextref#1{
      \global\advance\refno by 1
      \xdef#1{\the\refno}}

\def\bref {\ref\global\advance\refno by 1\key{\the\refno}}

%%%%%%%%% bibliographie  %%%%%%%%%%%

\nextref\BRU
\nextref\DEM 
\nextref\DUJ
\nextref\MCQ
\nextref\MOR
\nextref\NEV
\nextref\PAU
\nextref\SIK

%%%%%%%%% en tete%%%%%%%
\topmatter

 %\abstract
 %\endabstract
\title 
Singularit\'es des courants d'Ahlfors
 \endtitle

\author  Julien Duval
%\quad
%\the\day\ 
%\the\month\  
%\the\year
%\newline preliminary version, do not distribute
\endauthor

\abstract\nofrills{\smc R\'esum\'e. }\ \ On montre qu'une courbe alg\'ebrique charg\'ee par un courant issu d'une courbe enti\`ere
 est rationnelle ou elliptique. Ceci r\'epond \`a une question de M. Paun. 

\null
{\smc \hskip 2cm Singularities of Ahlfors currents }

\null
\noindent
{\smc Abstract.}\ We prove that an algebraic curve charged by a current coming from an entire curve
is rational or elliptic. This answers a question by M. Paun.
 \endabstract

\endtopmatter 

\document

\subhead 0. Introduction \endsubhead
\stno=1

\null
\footnote""{Mots-cl\'es : Kobayashi  hyperbolicity, entire curves, currents.}
\footnote""{Class.AMS : 32Q45, 32U40.}

Soit $X$ une vari\'et\'e compacte complexe hermitienne non hyperbolique au sens de Kobayashi. Elle poss\`ede
alors une courbe enti\`ere, i.e. une application holomorphe non constante de $\C$ dans $X$.

 Depuis Ahlfors (voir par exemple [\BRU]),
on sait associer \`a cette courbe enti\`ere un courant positif ferm\'e. Celui-ci s'obtient comme limite
de courants d'int\'egration normalis\'es $\frac{[\Dn]}{\text{aire}(\Dn)}$ sur des disques concentriques
de la courbe enti\`ere, choisis de sorte que la longueur de leur bord devienne n\'egligeable devant leur aire
: $\text{long}(\partial \Dn)=\text{o(aire}(\Dn))$.

Ces courants apparaissent notamment dans la preuve par M. McQuillan de la conjecture de Green-Griffiths pour
 certaines surfaces ([\MCQ], voir aussi [\BRU]). Celle-ci stipule qu'une courbe enti\`ere dans une vari\'et\'e alg\'ebrique
de type g\'en\'eral d\'eg\'en\`ere alg\'ebriquement.

\null

L'objet de cet article est de comprendre la partie singuli\`ere de ces courants.
Pr\'ecis\'ement, depuis Siu (voir par exemple [\DEM]), on sait d\'ecomposer tout courant positif ferm\'e en
 une partie singuli\`ere et une partie diffuse. La partie diffuse est peu singuli\`ere : elle a des nombres de Lelong nuls hors
 d'un ensemble d\'enombrable de points. La partie singuli\`ere rend compte des courbes analytiques charg\'ees par le courant :
 c'est une combinaison des courants d'int\'egration sur ces courbes. 

\null

A la suite de M. Paun (voir [\PAU] qui contient aussi des r\'eponses partielles), il est naturel de conjecturer qu'une courbe
 analytique charg\'ee par un courant
issu d'une courbe enti\`ere est rationnelle ou elliptique. En effet, la perte d'hyperbolicit\'e au sens de Kobayashi
due \`a la pr\'esence de la courbe enti\`ere doit se refl\'eter \`a l'endroit o\`u le courant se concentre le plus. Nous v\'erifions
 cette conjecture dans le cas o\`u la vari\'et\'e $X$ est alg\'ebrique. La m\'ethode repose sur une propri\'et\'e de
rel\`evement des arcs de la courbe alg\'ebrique \`a la courbe enti\`ere. Celle-ci est li\'ee au caract\`ere laminaire
du courant issu de la courbe enti\`ere (comparer avec [\DUJ]).

 Techniquement, il est commode d'\'elargir la classe des courants
consid\'er\'es en oubliant les courbes enti\`eres et en ne retenant que les suites
de disques, voire d'unions finies de disques, dont ils sont limites :

\null

\bf D\'efinition. \rm Un courant $T$ est dit {\it d'Ahlfors} s'il existe une suite $(\Dn)$ de r\'eunions finies de disques holomorphes dans $X$
telle que long($\partial \Dn)=o(\text{aire}(\Dn))$ et
 $T=\lim \frac{[\Dn]}{\text{aire}(\Dn)}$.
 
\null

Cette d\'efinition couvre aussi les courants issus de courbes enti\`eres apr\`es r\'egulari-sation \`a la Nevanlinna (voir [\MCQ] et
 [\BRU]).
Voici notre r\'esultat :

\thm{Th\'eor\`eme} Soit $T$ un courant d'Ahlfors dans une vari\'et\'e alg\'ebrique $X$. On suppose que
$T$ charge une courbe irr\'eductible $C$. Alors $C$ est rationnelle ou elliptique. \fthm

La d\'emonstration proc\`ede par l'absurde en supposant $C$ de genre au moins deux. On trouve alors un lacet $\gamma$ avec
un unique point double dans $C$ -une fi-gure huit- qui engendre dans le groupe fondamental de $C$ un groupe libre \`a
 deux g\'en\'erateurs. Comme $T$ se concentre sur $C$, on peut relever $\gamma$ en un graphe $\gamma_n$ dans $\Dn$ convergeant
 vers
$\gamma$ et poss\'edant beaucoup de cycles (de l'ordre de l'aire de $\Dn$). Comme $\Dn$ est simplement connexe,
 le graphe $\gamma_n$ va y d\'ecouper beaucoup de disques (de l'ordre de l'aire de $\Dn$).
 On exhibe ainsi une suite de disques holomorphes d'aire born\'ee a priori dont le bord converge vers un
 lacet non trivial de $\gamma$. Par un r\'esultat classique de compacit\'e, cette suite produit un disque limite dans
 $C$
bordant ce lacet, c'est la contradiction.

La propri\'et\'e de rel\`evement des arcs de $C$ \`a $\Dn$, cruciale, s'obtient en analysant $\Dn$ sous une projection
 centrale via la th\'eorie d'Ahlfors de recouvrement des surfaces (voir par exemple [\NEV]).
 C'est l\`a qu'intervient l'hypoth\`ese d'alg\'ebricit\'e de $X$.

\null

Avant de pr\'eciser ce sch\'ema, je tiens \`a remercier chaleureusement M. Paun pour m'avoir soumis ce probl\`eme.

\null

  \subhead 1. Pr\'eliminaires \endsubhead

\null
  
Enon\c cons pour commencer le peu de th\'eorie d'Ahlfors et le r\'esultat de compacit\'e des disques
 dont nous aurons besoin. Dans la suite, les espaces projectifs seront munis de leur m\'etrique standard.

\null

  a) {\bf Un peu de th\'eorie d'Ahlfors.}  
  
    \null
   
Nous renvoyons \`a la monographie [\NEV] pour les d\'etails de cette th\'eorie. Elle a \'et\'e cr\'e\'ee par Ahlfors pour
 traduire en termes g\'eom\'etriques la th\'eorie de Nevanlinna de distribution des valeurs.

\null
  
Pr\'ecis\'ement, soit $f:\Si \rightarrow \Si_0$ une application holomorphe entre deux surfaces de Riemann compactes
 \`a bord. On munit $\Si_0$ d'une m\'etrique qui induit via $f$ une pseudom\'etrique sur $\Si_0$. Le r\'esultat central
 de la th\'eorie est une in\'egalit\'e de Riemann-Hurwitz approch\'ee pour $f$ 
$$          \chi^-(\Si) \leq s\chi (\Si_0)+c_0l. $$

Ici, la caract\'eristique d'Euler n\'egative n\'eglige les composantes de $\Si$
qui sont des disques ou des sph\`eres,
$s$ est le nombre moyen de feuillets de $f$ : $s=\frac {\text{aire}(\Si)}{\text{aire}(\Si_0)}$,
 $l$ est la longueur du bord relatif de $f$ : $l=\text{long}(\partial \Si \setminus f^ {-1}(\partial \Si_0))$,
et la constante $c_0$ ne d\'epend que de $\Si_0$ et de sa m\'etrique.

\null

Ahlfors en d\'eduit son th\'eor\`eme des \^\i les, l'analogue du second th\'eor\`eme fondamental de Nevanlinna :

\thm{Th\'eor\`eme des \^\i les} Soit $f:\De \rightarrow \p$ une application holomorphe entre une r\'eunion finie $\De$
 de disques et
 $\p$.
 On se fixe $k$  disques topologiques
disjoints dans $\p$ et on appelle  \^\i les  de $f$ les composantes connexes de la pr\'eimage des disques sur lesquelles $f$
 est propre.
Alors le nombre $i$ d'\^\i les v\'erifie $i \geq s(k-2)-cl$.

\fthm

Ici, comme plus haut, $s$ d\'esigne le nombre de feuillets moyen de $f$, $l$ la longueur du bord de $\De$, et 
$c$ ne d\'epend que de la configuration de disques topologiques.

Il suffit en effet d'appliquer l'in\'egalit\'e de Riemann-Hurwitz \`a  $\Si_0=\p$ priv\'e des k disques et $\Si =f^{-1}(\Si_0)$,
 sachant
 que le nombre de feuillets moyen de $f$ sur $\Si_0$ diff\`ere de $s$ par un terme contr\^ol\'e par $l$.

\null

Ceci entra\^\i ne qu'une suite de fonctions m\'eromorphes sur des disques 
se comporte asymptotiquement comme un rev\^etement si la longueur du bord devient n\'egligeable devant l'aire.
 Elle poss\`ede la propri\'et\'e de rel\`evement suivante :
\thm{Rel\`evement} Soit $f_n:\Dn \rightarrow \p$ une suite d'applications holomorphes entre des r\'eunions finies $\Dn$ de disques et $\p$.
 On suppose, avec les notations pr\'ec\'edentes, que $l_n=o(s_n)$. On se fixe $\epsilon >0$ et des
 disques topologiques disjoints dans $\p$ en nombre
 $k>4/\epsilon$. Alors un de ces disques topologiques
poss\`ede au moins $s_n(1-\epsilon)$ rel\`evements via $f_n$ pour $n$ assez grand.
\fthm

Ici, un rel\`evement est une composante connexe de la pr\'eimage du disque sur laquelle $f_n$ est bijective.

 En effet, par le th\'eor\`eme
 des \^\i les, un des disques topologiques va poss\'eder au moins $s_n(1-\frac{\epsilon}{2})$ \^\i les dans sa pr\'eimage pour $n$
 assez grand. Or la restriction
 de $f_n$ \`a chaque
 \^\i le est un rev\^etement ramifi\'e au-dessus du disque. Ce rev\^etement sera donc de degr\'e 1 pour au moins
 $s_n(1-\epsilon)$ \^\i les, puisque le nombre moyen de feuillets de $f_n$ sur ce disque est asymptotiquement $s_n$.
  
   \null

  b) {\bf Compacit\'e de disques holomorphes.}
 
 \null 

Voici l'\'enonc\'e qui nous sera utile.
Il remonte au moins \`a Gromov dans un contexte bien plus large (voir par exemple [\SIK]).
\thm {Compacit\'e}
Soit $f_n:D \rightarrow \pp$ une suite d'applications holomorphes du disque unit\'e $D$ dans $\pp$.
 On suppose l'aire de $f_n(D)$ uniform\'ement born\'ee. Alors, quitte \`a extraire, 

 i) $f_n$ converge localement uniform\'ement vers $f_\infty$ sur $D$ priv\'e d'un nombre fini
de points d'explosion;

ii) $f_\infty$ se prolonge en un disque holomorphe encore not\'e $f_\infty:D\rightarrow \pp$;

iii) en un point d'explosion $e$, il existe une suite de disques $d_n$ dans $D$ tendant vers $e$ telle que
$f_n(d_n)$ converge au sens de Hausdorff vers une courbe rationnelle (une bulle).

\fthm

En particulier, si l'aire de $f_n(D)$ reste born\' ee par $\frac{1}{2}$ il ne se produira pas d'explosion puisqu'une
courbe rationnelle est d'aire au moins 1. On aura dans ce cas convergence uniforme locale de $f_n$ vers $f_\infty$ sur $D$ apr\`es extraction.

\null

 \subhead 2. D\'emonstration du th\'eor\`eme \endsubhead
 
 \null
 
Notons que, par plongement de $X$ dans un espace projectif, il suffit de regarder les courants d'Ahlfors
de $P^N(\C)$ chargeant une courbe projective. La d\'emonstration proc\`ede en trois \'etapes : la premi\`ere
est une r\'eduction \`a $\pp$, la seconde traite de la propri\'et\'e de rel\`evement des arcs, tandis que la
troisi\`eme d\'etaille la fin de l'argument.

 Pour des pr\'ecisions sur les courants et leurs nombres de Lelong, nous renvoyons
 le lecteur
\`a [\DEM]. Dans la suite on ne distinguera plus les disques param\'etr\'es de leur image g\'eom\'etrique. En particulier
les aires des images seront compt\'ees avec multiplicit\'e.
 Que le lecteur nous pardonne cette libert\'e.

\null

a) \bf R\'eduction de dimension.\rm

 \null

Soit donc un courant d'Ahlfors $T$ de $P^N(\C)$ chargeant une courbe projective $C$.
On projette $T$ dans un hyperplan projectif par une projection $\pi$ de centre $p$. Si $p$ est choisi assez g\'en\'erique,
 le nombre de Lelong de $T$ en $p$ sera nul et la courbe $\pi(C)$ aura m\^eme genre que $C$ tant que $N \geq 3$.
Le lemme suivant permet donc de se ramener au cas $N=2$ par projections centrales successives.

\thm{Lemme} Le courant projet\'e $\Ti$ est un courant d'Ahlfors chargeant la courbe projective $\pi(C)$.
\fthm

 \demo {D\'emonstration} La difficult\'e r\'eside dans la singularit\'e de $\pi$ en $p$. Pr\'ecisons d'abord la d\'efinition de $\Ti$ : 
si $B_r$ d\'esigne la boule de centre $p$ et de rayon $r$, on pose $\Ti=\lim_{r\rightarrow 0}\pi_*(T\vert_{(B_r)^c})$.
 Clairement $\Ti$
charge $\pi(C)$. Reste \`a v\'erifier que c'est un courant d'Ahlfors.
 
Pour cela, partons de la suite de r\'eunions finies $\Dn$ de disques holomorphes donnant $T$. Elle satisfait 
$\lim \frac{[\Dn]}{a_n}=T$ et
 long($\partial \Dn)=
o(a_n)$ avec  $a_n=\text{aire}(\Delta_n)$. La suite analogue $\DTn$ pour $\Ti$ s'obtient essentiellement en projetant par $\pi$
 la partie de $\Dn$ loin de $p$.
 Le contr\^ole de la longueur de son bord va provenir de l'annulation du nombre de Lelong de $T$ en $p$. Pr\'ecisons ceci.

Soit $r>0$ fix\'e provisoirement. Dans ce qui suit, les in\'egalit\'es seront \`a constante multiplicative ind\'ependante de $r$ et $n$ pr\`es.
Voici comment se traduit la nullit\'e du nombre de Lelong de $T$ en $p$ pour $n$ assez grand :
$$ \text{aire}(\Dn \cap B_r)\leq\epsilon(r)r^2a_n.\tag1 $$
Ici $\epsilon(r)$ tendra vers $0$ avec $r$.

L'in\'egalit\'e de coaire (voir par exemple [\MOR]) fournit alors $r_n$ avec $\frac{r}{2}<r_n<r$ et 
$$ \text{long}(\Dn \cap \partial B_{r_n}) \leq \epsilon(r)ra_n.\tag2$$

Remarquons que $\Dn \cap (B_{r_n})^c$ est une r\'eunion de disques trou\'es. On rebouche ces trous par les composantes
connexes correspondantes de $\Dn \cap B_{r_n}$. On obtient ainsi une r\'eunion $\Dn'$ de disques dont le bord consiste
en deux parties, l'une incluse dans $\Dn \cap \partial B_{r_n}$, l'autre \'etant $\partial \Dn \cap  (B_{r_n})^c$.
Notons $\DTn$ la projection de $\Dn'$ par $\pi$. C'est la suite recherch\'ee.

 Puisque $\pi$ dilate au plus les longueurs d'un facteur $\frac{1}{r}$ sur $(B_{r_n})^c$, on a, d'apr\`es (2),
$$ \text{long}(\partial \DTn)\leq \epsilon(r)a_n \tag3$$
pour $n$ assez grand.

Comparons maintenant $a_n$ avec $\tilde a_n$ l'aire de $\DTn$. Pour cela, notons $\omega$ et $\widetilde{\omega}$ les formes de Fubini-Study
 respectives
de $P^N(\C)$ et de l'hyperplan sur lequel on projette. On a $a_n=\int_{\Dn}\omega$ et $\tilde a_n=\int_{\Dn'}\pi^*\widetilde{\omega}$.
  Remarquons que $\pi^*\widetilde{\omega}-\omega=dd^cu$ o\`u $u$ est une fonction lisse hors de $p$
avec une singularit\'e logarithmique en $p$. En particulier $d^cu$ est de l'ordre de $\frac{1}{r}$ sur $(B_{r_n})^c$.
Il s'ensuit, par le th\'eor\`eme de Stokes et les estim\'ees (1) et (3), que
$$\vert \tilde a_n-a_n\vert \leq \int_{\Dn \setminus \Dn'} \omega+\vert \int_{\partial \Dn'}d^cu\vert \leq
\text{aire}(\Dn \cap B_{r_n})+\text{long}(\partial \DTn)\leq \epsilon(r)a_n$$
 pour $n$ assez grand. Donc long($\partial \DTn)\leq \epsilon(r)\tilde a_n$.

Faisons tendre maintenant $r$ vers 0. Apr\`es extraction diagonale sur les suites $(\DTn)$ produites \`a $r$ fix\'e,
on v\'erifie que $\Ti=\lim \frac{[\DTn]}{\tilde a_n}$ et
 long($\partial \DTn)=
o(\tilde a_n)$. C'est donc bien un courant d'Ahlfors. $\square$

\enddemo
\null
\noindent
\bf Remarque. \rm Le
m\^eme type d'estim\'ee montre que l'aire de $\pi(\Dn' \cap B_{r_n})$ est contr\^ol\'ee par $\epsilon(r)a_n$ pour $n$ assez grand.

\null

b) \bf Rel\`evement des arcs. \rm

 \null

 On s'est ramen\'e \`a un courant d'Ahlfors $T$ de $\pp$ chargeant une courbe projective $C$.
Pr\'ecis\'ement $T\vert_C=\nu[C]$ avec $\nu >0$.

 Analysons encore
 $T$ sous une projection $\pi$ de centre $p$ sur une droite projective $L$. Si $p$ est choisi assez g\'en\'erique,
on modifie comme plus haut les r\'eunions $\Dn$ de disques d\'efinissant $T$ pour que les longueurs des
bords de $\Dn$ et $\pi(\Dn)$ deviennent n\'egligeables devant l'aire de $\Dn$, celle-ci soit \'equivalente \`a
l'aire de $\pi(\Dn)$, et $$\lim_r \limsup_n \frac{\text{aire}(\pi(\Dn \cap B_r))}{\text{aire}(\Dn)}=0. \tag4$$

Notons $s_n$ l'aire de $\pi(\Dn)$. C'est le nombre moyen de feuillets de $\pi:\Dn \rightarrow L$.
Dans la suite on appellera {\it arc} un chemin injectif dans $C$ qui se projette injectivement par $\pi$ dans $L$.
On dira qu'un arc $\alpha$ est $\epsilon$-{\it relevable} si on peut l'\'epaissir en une bande $\beta$ dans $C$
poss\'edant au moins $s_n(\nu-\epsilon)$ relev\'es dans $\Dn$ apr\`es extraction, qui convergent vers $\beta$. On parlera
de relev\'es {\it  proches}.
Ici un relev\'e de $\beta$ est une composante $\beta'$ de $\pi^{-1}(\pi(\beta))$ dans $\Dn$ telle que $\pi: \beta' \rightarrow
\pi(\beta)$ soit bijective.

\thm{Lemme} Tout arc $\alpha$ de $C$ peut se perturber en un arc $\epsilon$-relevable.
\fthm

{\it D\'emonstration} (comparer \`a [\DUJ]). On se fixe un voisinage $V$ de $\alpha$ dans $C$ qui se projette injectivement par $\pi$ sur $U$.
 S'il est suffisamment
fin, la masse de $T$ dans $\pi^{-1}(U)$ est inf\'erieure \`a $\frac{\epsilon}{4}$.
 Soient maintenant $k$ perturbations disjointes de $\alpha$ dans $V$, avec
$k> \frac{8}{\epsilon}$. On les \'epaissit un peu en $k$ fines bandes disjointes.

 Appliquons la propri\'et\'e
 de rel\`evement 1.a) \`a $\pi : \Dn \rightarrow L$ et aux $k$ disques disjoints qui sont les images par $\pi$
de ces bandes. On obtient qu'une de ces bandes $\beta_0$ a une projection $\pi(\beta_0)= \tilde \beta$ poss\'edant $s_n(1-\frac{\epsilon}{2})$
relev\'es dans $\Dn$ pour $\pi$, si $n$ est assez grand et apr\`es extraction. Parmi ceux-ci, au moins $s_n(1-\epsilon)$
sont d'aire inf\'erieure \`a $1/2$ par le contr\^ole de la masse de $T$ dans $\pi^{-1}(U)$. On les note $\beta_{n,i}$.
 Reste \`a en localiser une partie pr\`es de $\beta_0$.

Introduisons pour cela un peu de terminologie.
Soient $Q$ le c\^one $\pi^{-1}(\tilde \beta)$ et
 $B$ l'ensemble des disques holomorphes $\beta$ de $Q$ d'aire inf\'erieure \`a
 1/2 et se projetant bijectivement par $\pi$ sur $\tilde \beta$. On peut voir les disques de $B$ comme des sections de $\pi$ au-dessus de
 $\tilde \beta$.

Par 1.b) une suite dans $B$ poss\`edera toujours une suite extraite convergente puisque la borne d'aire interdit les explosions.
 La limite est encore dans $B$ ou vaut $p$.
Ainsi
$B$ est relativement compact pour la convergence uniforme locale sur $\tilde \beta$, et l'espace $M$ des mesures positives sur $B$
 de masse au plus 1
est compact pour la convergence faible.

A une mesure $\lambda$ de $M$ on associe un courant ``g\'eom\'etrique'' $T_{\lambda}$ dans $Q$ par la formule
 $T_{\lambda}=\int_{B} [\beta]d\lambda(\beta)$. Cette application de $M$ sur l'espace $G$ des courants g\'eom\'etriques \'etant continue,
 on en d\'eduit
 que
 $G$ est compact pour la topologie faible.

 Revenons maintenant \`a nos relev\'es $\beta_{n,i}$.
Par construction, la suite
 $(\frac{1}{s_n} \sum_i[\beta_{n,i}])$ est une suite
 de courants g\'eom\'etriques $(T_{\lambda_n})$ pour $\lambda_n = \frac{1}{s_n} \sum_i\delta_{\beta_{n,i}}$.
Elle converge apr\`es extraction vers $T_\lambda$ dans $B$.
Comme $T_{\lambda_n} \leq \frac{1}{s_n}[\Dn \cap Q]$ on a $T_{\lambda}\leq T\vert_Q$.

 Par ailleurs $\lambda$ a une masse minor\'ee
par $1-\epsilon$ car le nombre de disques $\beta_{n,i}$ s'\'echappant vers $p$ est n\'egligeable d'apr\`es (4).
 Donc $\pi_*(T_{\lambda})> (1-\epsilon)[\tilde \beta].$
 Comme $\pi_*(T\vert_Q)=[\tilde \beta]$ on en d\'eduit $$\pi_*(T\vert_Q-T_\lambda)<\epsilon[\tilde 
\beta].$$
Par hypoth\`ese $T$ charge $C$ avec un poids $\nu$ donc $T\vert_{\beta_0}= \nu[\beta_0]$. Par construction $T_\lambda\vert_{\beta_0}=\mu[\beta_0]$
 o\`u $\mu$ est la charge de
$\beta_0$ pour $\lambda$.
 Donc $$\pi_*(T\vert_Q-T_\lambda)\geq(\nu-\mu)[\tilde \beta].$$ Ainsi la charge de $\beta_0$ pour $\lambda$ est sup\'erieure \`a
 $\nu-\epsilon$. Il en est de m\^eme de la charge
pour $\lambda_n$ de tout voisinage de $\beta_0$ dans $B$ pour $n$ assez grand.
Autrement dit, au moins  $s_n(\nu-\epsilon)$ disques $\beta_{n,i}$
convergent vers $\beta_0$.  $\square$

\null

Un {\it graphe} $\gamma$ de $C$ est une r\'eunion d'un nombre fini $a$ d'arcs $\alpha_j$ se coupant en un nombre fini $i$ de points.
 Deux arcs sont disjoints ou se coupent exactement une fois, soit en leurs extr\'emit\'es (arcs cons\'ecutifs),
 soit en leurs int\'erieurs (arcs en croix).
 Le graphe $\gamma$ sera dit
 {\it relevable} si ses arcs sont $\epsilon$-relevables pour $\epsilon=\frac{\nu}{2a+6i}$. Gr\^ace au lemme on peut toujours perturber un graphe
en un graphe relevable en pr\'eservant sa combinatoire d'incidence. 
Notons $\Gamma$ la r\'eunion des bandes $\beta_j$ \'epaississant les arcs et $\Gamma_n$ la r\'eunion
de leurs relev\'es proches dans $\Dn$. C'est l'\'epaississement d'un graphe $\gamma_n$ dans $\Dn$ qui converge vers $\gamma$.

\thm{Fait} On a $\chi(\Gamma_n) \leq s_n(\chi(\Gamma)+\frac{1}{2})\nu$ pour $n$ assez grand.
\fthm

En effet les bandes $\beta_j$ sont disjointes ou se coupent sur des disques qui \'epais-sissent les intersections des $\alpha_j$. Donc
$\chi(\Gamma)=a-i$.
Par ailleurs les relev\'es proches des bandes sont disjoints ou se recollent sur des relev\'es proches de ces disques. Donc
$\chi(\Gamma_n)=a_n-i_n$ o\`u $a_n$ d\'esigne le nombre total de relev\'es proches de bandes et $i_n$ le nombre total de recollements.

Estimons ces quantit\'es. Pour cela remarquons que les relev\'es proches d'une bande $\beta_j$  donn\'ee (ou d'un disque) sont
 en nombre major\'e par
$s_n(\nu+\epsilon)$ pour $n$ assez grand. On aurait sinon $T\vert_{\beta_j} \geq (\nu+\epsilon)[\beta_j]$ par passage \`a la limite.

Donc $a_n \leq s_n(\nu+\epsilon)a$. Par ailleurs, au-dessus d'un disque intersection de deux bandes, on a au moins $s_n(\nu-\epsilon)$ relev\'es
proches provenant de chacune des bandes. Il s'ensuit n\'ecessairement au moins $s_n(\nu-3\epsilon)$ recollements au-dessus de ce disque.
Ainsi $i_n \geq s_n(\nu-3\epsilon)i$. On obtient donc par notre choix de $\epsilon$
  $$\chi(\Gamma_n)\leq s_n( \chi(\Gamma)\nu + \epsilon(a+3i))\leq s_n(\chi(\Gamma)+\frac{1}{2})\nu. $$

c) \bf Fin de l'argument. \rm

\null

Supposons par l'absurde $C$ de genre au moins 2. Soit $\gamma$ une figure huit (un lacet avec un unique point double) dans $C$
dont le groupe fondamental s'injecte dans celui de $C$. On reprend ici la terminologie du paragraphe pr\'ec\'edent.
 On peut toujours supposer apr\`es perturbation que $\gamma$ est un graphe
relevable.

Comme $\chi(\gamma)=-1$ on a, par l'estim\'ee pr\'ec\'edente, $\chi(\gamma_n)\leq -\frac{s_n\nu}{2}$ pour $n$ assez
grand.
Par un argument de dualit\'e d'Alexander, $\gamma_n$ d\'ecoupe donc dans $\Dn$ au moins $\frac{s_n\nu}{2}$ disques.
Or $s_n$ \'equivaut \`a l'aire de $\Dn$. En s\'electionnant pour chaque $n$ le disque $\dn$
composante de $\Dn \setminus \gamma_n$ d'aire minimale, on produit ainsi une suite de disques holomorphes d'aire born\'ee a priori.

Analysons la convergence de $(\dn)$. Commen\c cons par leur bord $\partial \dn$. Il est constitu\'e d'un cycle
de relev\'es proches d'arcs de $\gamma$. Par construction ceux-ci s'\'epaississent dans $\Dn$ en relev\'es proches des demi-bandes
correspondantes dans $\Gamma$.
 Ces relev\'es forment un anneau $A_n$ dans $\dn$ bordant $\partial \dn$ et d'aire comparable \`a la
longueur de $\partial \dn$. Il s'ensuit que le nombre de relev\'es proches d'arcs de $\gamma$ dans $\partial \dn$ est born\'e a priori. Donc,
 quitte \`a extraire, $\partial \dn$ converge vers un cycle non trivial d'arcs de $\gamma$. De m\^eme $A_n$ converge vers un anneau $A$
 immerg\'e dans $C$.
 En particulier le module de $A_n$ reste minor\'e.
Cela permet de choisir un param\'etrage holomorphe $f_n$ des disques $\dn$ par le disque unit\'e $D$ tel que la pr\'eimage de $A_n$ contienne
un anneau fixe $(r<\vert z\vert <1)$.

 Appliquons l'\'enonc\'e de compacit\'e 1.b). Comme $C$ n'est pas rationnelle, les \'eventuelles
explosions ne peuvent se produire que dans $(\vert z \vert \leq r)$. Ainsi le disque limite $\delta_\infty=f_\infty (D)$ contient $A$.
 Par prolongement analytique il est enti\`erement dans $C$.
 Or $\partial \delta_\infty$ est un lacet non trivial de $\gamma$ donc non homotope \`a z\'ero dans $C$.
C'est la contradiction.  $\square$

 \Refs

%\medskip
\widestnumber\no{99}
\refno=0

\bref \by M. Brunella \paper Courbes enti\`eres et feuilletages holomorphes \jour
Ens. Math. \vol45\yr1999\pages195--216
\endref

\bref \by J.-P. Demailly \book Monge-Amp\`ere operators, Lelong numbers and intersection theory,  {\rm in}
 Complex analysis and geometry \pages115--193 \bookinfo Univ. Ser. Math.\publ Plenum \yr1993 \publaddr New-York
\endref

\bref \by R. Dujardin \paper Laminar currents and dynamics I : Structure properties of laminar currents on $ P^2$
\jour preprint 2004 arXiv math.CV/0403292
\endref

\bref \by M. McQuillan \paper Diophantine approximations and foliations \jour Publ. Math. IHES \vol87\yr1998\pages121--174
\endref

\bref \by F. Morgan \book Geometric measure theory. A beginner's guide \publ Academic Press \yr 2000 \publaddr San Diego
\endref

\bref \by R. Nevanlinna \book Analytic functions \bookinfo Grund. Math. Wiss. \vol162 \publ Springer \yr 1970 \publaddr Berlin
\endref

\bref \by M. Paun \paper Currents associated to transcendental entire curves on compact K\"ahler manifolds \jour preprint 2003
\endref

\bref \by J.-C. Sikorav \book Some properties of holomorphic curves in almost complex manifolds, {\rm in} Holomorphic curves in symplectic geometry
\pages165--189 \bookinfo Prog. Math. \vol117 \publ Birkh\"auser \yr1994 \publaddr Basel
\endref

\endRefs

\address 
\noindent  
Laboratoire \'Emile Picard, 
  Universit\'e Paul Sabatier, 31062 Toulouse Cedex 09.
 \endaddress
\email 
  duval\@picard.ups-tlse.fr 
\endemail

\enddocument